\documentclass[11pt]{article}
\usepackage{amssymb, amsthm, amsmath, amscd}
\newtheorem{thm}{Theorem}[section]
\newtheorem{lem}[thm]{Lemma}
\newtheorem{prop}[thm]{Proposition}
\newtheorem{cor}[thm]{Corollary}
\theoremstyle{definition}\newtheorem{df}[thm]{Definition}
\theoremstyle{definition}\newtheorem{rem}[thm]{Remark}
\theoremstyle{definition}

\renewcommand{\phi}{\varphi}

\newcommand{\Z}{\mathbb{Z}}
\newcommand{\Q}{\mathbb{Q}}

\newcommand{\T}{\mathbb{T}}

\newcommand{\Aut}{\operatorname{Aut}}

\newcommand{\hm}{homomorphism}
\newcommand{\af}{\alpha}

\newcommand{\dt}{\delta}
\newcommand{\ep}{\epsilon}
\newcommand{\andeqn}{\,\,\,{\rm and}\,\,\,}
\newcommand{\rforal}{\,\,\,{\rm for\,\,\,all}\,\,\,}
\newcommand{\CA}{$C^*$-algebra}
\newcommand{\SCA}{$C^*$-subalgebra}
\newcommand{\tr}{{\rm TR}}
\newcommand{\beq}{\begin{eqnarray}}
\newcommand{\eneq}{\end{eqnarray}}
\newcommand{\tforal}{\,\,\,\text{for\,\,\,all}\,\,\,}

\title{ The Rokhlin Property for Automorphisms on
Simple \CA s}
\author{Huaxin Lin }
\date{}

\begin{document}

\maketitle

\begin{abstract}
We study  a general Kishimoto's problem for automorphisms on
simple \CA s with tracial rank zero. Let $A$ be a unital separable
simple \CA\, with tracial rank zero and let $\alpha$ be an
automorphism. Under the assumption that $\alpha$ has certain
Rokhlin property, we present a proof that $A\rtimes_{\alpha}\Z$
has tracial rank zero. We also show that if the induced map
$\alpha_{*0}$ on $K_0(A)$ fixes a ``dense" subgroup of $K_0(A)$
then the tracial Rokhlin property implies a stronger Rokhlin
property. Consequently, the induced crossed product \CA s have
tracial rank zero.
\end{abstract}

\section{Introduction}

The Rokhlin property in ergodic theory was adopted to the context
of von Neumann algebras by Connes (\cite{Con}). The Rokhlin
property (with various versions) was also introduced to the study
of automorphisms on \CA s (see, for example, Herman and Ocneanu
(\cite{HO}), R\o rdam (\cite{R}, Kishimoto  (\cite{K1}) and
Phillips \cite{P1} among others-- see also the next section).

A conjecture of Kishimoto can be formulated as follows: Let $A$ be
a unital simple A$\T$-algebra of real rank zero and $\alpha$ be an
approximately inner automorphism. Suppose that $\alpha$ is
``sufficiently outer", then the crossed product of the
A$\T$-algebra by $\alpha,$ $A\rtimes_{\alpha}\Z,$ is again a
unital A$\T$-algebra. Kishimoto showed that this is true for a
number of cases, in particular, for some cases that  $A$ has a unique
tracial state.

   Kishimoto proposed that the appropriate notion of outerness is
   the Rokhlin property (\cite{K2}).
Kishimoto's problem has a more general setting:

{\bf P1}\,\,\,  Let $A$ be a unital separable simple
\CA\, with tracial rank zero and let $\alpha$ be an automorphism.
Suppose that $\alpha$ has a Rokhlin property. Does
$A\rtimes_{\alpha}\Z$ have tracial rank zero?

 If in addition $A$ is assumed to be amenable and satisfy the Universal
 Coefficient Theorem, then, by the classification theorem
 (\cite{Lnduke}), $A$ is an AH-algebra with slow dimension growth and
 with real rank zero.
 If $A\rtimes_{\alpha}\Z$ has tracial rank zero, then, again,
 by \cite{Lnduke}, $A\rtimes_{\alpha}\Z$ is an AH-algebra
 (with slow dimension growth and with real rank zero).
Note that  simple A$\T$-algebras with real rank zero are exactly
 those simple AH-algebras with torsion free $K$-theory, with
 slow dimension growth and with real rank zero.
 It should also be noted (\cite{Lnpams}) that  a unital simple AH-algebra
 has slow dimension growth and real rank zero if and only if
 it has tracial rank zero.

If $K_i(A)$ are torsion free and $\alpha$ is approximately inner,
then $K_i(A\rtimes_{\alpha}\Z)$ is torsion free. Thus an
affirmative answer to the problem {\bf P1} proves the original
Kishimoto's conjecture (see also \cite{LO}). One should also
notice that if $\alpha$ is not approximately inner,
$K_i(A\rtimes_{\alpha}\Z)$ may have torsion even if $A$ does not.
Therefore, in Kishimoto's problem, the restriction that $\alpha$
is approximately inner can not be removed. So it is appropriate to
replace $A\T$-algebras by AH-algebras if the requirement that
$\alpha$ is approximately inner is removed.

In this paper we report some of the recent development on this
subject.

Let $A$ be a unital separable simple \CA\, with tracial rank zero
and let $\alpha$ be an automorphism on $A.$ In section 3, we
present a proof that if $\alpha$ satisfies the tracial cyclic
Rokhlin property (see \ref{ID2} below) then the crossed product
$A\rtimes_{\alpha}\Z$ has tracial rank zero which gives a solution
to {\bf P1}, provided that $[\alpha^r]=[{\rm id}_A]$ in $KL(A,A).$
  In section 4, we discuss when $\alpha$ has the
tracial cyclic Rokhlin property. The tracial Rokhlin property
introduced by N. C. Phillips (see \cite{P1} and \cite{OP}) has
been proved to be a natural generalization of original Rokhlin
towers for ergodic actions. We prove that, if in addition,
$\alpha^r_{*0}|_{G}={\rm id}|_G$ (for some integer $r\ge 1$) for
some subgroup $G\subset K_0(A)$ for which $\rho_A(G)$ is dense in
$\rho_A(K_0(A)),$ then that $\alpha$ has tracial Rokhlin property
implies that $\alpha$ has tracial cyclic Rokhlin property.

\vspace{0.2in}

{\bf Acknowledgement} This paper was intended for the Proceedings
of GPOTS 2005. Most of it was written when the author was in East
China Normal University during the summer 2005. It is partially
supported by a grant from NSF of U.S.A, Shanghai Priority Academic
Disciplines and Zhi-Jiang Professorship from East China Normal
University.

\section{Preliminaries}

\vskip 3mm

We will use the following convention:

(i) Let $A$ be a \CA, let $a\in A$ be a positive element and let
$p\in A$ be a projection. We write $[p]\le [a]$ if there is a
projection $q\in {\overline{aAa}}$ and a partial isometry $v\in A$
such that $v^*v=p$ and $vv^*=q.$

\vspace{0.12in}

(ii) Let $A$ be a \CA. We denote by $Aut(A)$ the automorphism
group of $A.$ If $A$ is unital and $u\in A$ is a unitary, we
denote by ${\rm ad}\, u$ the inner automorphism defined by ${\rm
ad}\, u(a)=u^*au$ for all $a\in A.$ \vspace{0.12in}

(iii) Let $T(A)$ be the tracial state space of a unital \CA\, $A.$
It is a compact convex set. Denote by $Aff(T(A))$ the normed space of
all real affine continuous functions on $T(A).$ Denote by $\rho_A:
K_0(A)\to Aff(T(A))$ the \hm\, induced by
$\rho_A([p])(\tau)=\tau(p)$ for $\tau\in T(A).$

It should be noted, by \cite{BH}, if $A$ is a unital simple
amenable \CA\, with real rank zero and stable rank and weakly
unperforated $K_0(A),$ $\rho_A(K_0(A))$ is dense in $Aff(T(A)).$
\vspace{0.12in}

(iv)  Let $A$ and $B$ be two \CA s and $\phi, \psi: A\to B$ be two
maps. Let $\ep>0$ and ${\cal F}\subset A$ be a finite subset.
 We write
$$
\phi\approx_{\ep}\psi \,\,\,{\rm on}\,\,\, {\cal F},
$$
if
$$
\|\phi(a)-\psi(a)\|<\ep \rforal a\in {\cal F}.
$$

\vspace{0.12in}

(v) Let  $x\in A,$ $\ep>0$ and ${\cal F}\subset A.$ We write
$x\in_{\ep} {\cal F},$ if ${\rm dist}(x, {\cal F})<\ep,$ or there
is $y\in {\cal F}$ such that $\|x-y\|<\ep.$

\vspace{0.12in}

(vi) If $h: A\to B$ is a \hm, then $h_{*i}: K_i(A)\to K_i(B)$
($i=0,1$) is the induced \hm.

\vskip 3mm

We recall the definition of tracial topological rank of
C*-algebras.

\begin{df}\label{ID0}{\rm \cite[Theorem 6.13]{Lntr0}}
{\rm Let $A$ be a unital simple \CA~. Then $A$ is said to  have
{\it tracial (topological) rank zero} if  for any finite set
${\mathcal F} \subset A$, and $\ep > 0$ and any non-zero positive
element $a \in A$, there exists a finite dimensional \SCA\, $B
\subset A$ with  ${\rm id}_B = p$ such that
\begin{itemize}
\item[$(1)$] $\|px-xp\| < \ep$ for all $x \in {\mathcal F}$,
\item[$(2)$] $pxp \in_\ep B$ for all $x \in {\mathcal F}$,
\item[$(3)$] $[1 - p] \leq [a]$.
\end{itemize}
If $A$ has tracial rank zero, we write $\tr(A)=0$. }
\end{df}

If $A$ is assumed to have the Fundamental Comparison Property
(i.e., for any two projections $p,\,q\in A,$ $\tau(p)<\tau(q)$ for
all $\tau\in T(A)$ implies that $p$ is equivalent to a projection
$p'\le q.$), then the third condition may be replaced by
$\tau(1-p)<\ep$ for all $\tau\in T(A).$ It is proved
(\cite{Lntr0}, or see 3.7 of \cite{Lnbk}) that if $TR(A)=0,$ then
$A$ has the Fundamental Comparison Property, as well as real rank
zero and stable rank one. Every simple AH-algebra with real rank
zero and with the Fundamental Comparison Property has tracial rank
zero (\cite{Lnpams} and \cite{EG}). Other simple \CA s with
tracial rank zero may be found in \cite{Lncrel}.

There are several versions of the Rokhlin property (see \cite{HO},
\cite{R}, \cite{K1}  and \cite{Iz}).

The following is defined in \cite[Definition 2.1]{OP}.

\begin{df}\label{ID1}
{\rm Let $A$ be a simple unital \CA~ and let $\alpha \in \Aut(A)$.
We say $\alpha$ has the {\it tracial Rokhlin property} if for
every finite set ${\cal F} \subset A$, every $\ep > 0$, every $n
\in {\mathbb N}$, and every nonzero positive element $a \in A$,
there are mutually orthogonal projections $e_1,e_2, \dots, e_n
\in A$ such that:
\begin{itemize}
\item[$(1)$] $\| \alpha (e_j) - e_{j + 1} \| < \ep$ for $1 \leq j
\leq n - 1$. \item[$(2)$] $\| e_j a - a e_j \| < \ep$ for $0 \leq
j \leq n$ and all $a \in {\cal F}$. \item[$(3)$] With $e = \sum_{j
= 0}^{n} e_j$, $[1 - e] \leq [a]$.
\end{itemize}
}
\end{df}

The following result of Osaka and Phillips is the tracial Rokhlin
version of Kishimoto's result in the case of simple unital
$A{\mathbb T}$-algebras with a unique trace \cite[Theorem
2.1]{K2}.

\begin{thm}\label{IT1}{\rm cf. \cite{OP}}
{\rm Let $A$ be a simple unital \CA~ with $\tr(A) = 0,$ and
suppose that $A$ has a unique tracial state. Then the following
conditions are equivalent:
\begin{itemize}
\item[$(1)$] $\alpha$ has the tracial Rokhlin property.
\item[$(2)$] $\alpha^m$ is not weakly inner in the GNS
representation $\pi_\tau$ for any $m \not= 0$. \item[$(3)$] $A
\rtimes_\alpha{\mathbb Z}$ has real rank zero. \item[$(4)$] $A
\rtimes_\alpha{\mathbb Z}$ has a unique trace.
\end{itemize}
}
\end{thm}

\vskip 3mm

\vskip 3mm

We define a stronger version of the tracial Rokhlin property
similar to the approximately Rokhlin property in \cite[Definition
4.2]{K2}.

\vskip 3mm

\begin{df}\label{ID2}
{\rm Let $A$ be a simple unital
 \CA~ and
let $\alpha \in \Aut(A)$. We say $\alpha$ has the {\it tracial
cyclic Rokhlin property} if for every finite set ${\cal F} \subset
A$, every $\ep > 0$, every $n \in {\mathbb N}$, and every nonzero
positive element $a \in A$, there are mutually orthogonal
projections $e_0, e_1, \dots, e_n \in A$ such that
\begin{itemize}
\item[$(1)$] $\| \alpha (e_j) - e_{j + 1} \| < \ep$ for $0 \leq j
\leq n$, where $e_{n + 1} = e_0$. \item[$(2)$] $\| e_j a - a e_j
\| < \ep$ for $0 \leq j \leq n$ and all $a \in {\cal F}$.
\item[$(3)$] With $e = \sum_{j = 0}^{n} e_j$, $[1 - e] \leq [a]$.
\end{itemize}
}
\end{df}
The following is a restatement of Theorem 3.4 of \cite{Lncd}.

\begin{thm}\label{CDT}
Let $C$ be a unital AH-algebra  and let $A$ be a unital simple \CA\, with tracial rank zero.
Suppose that  $\phi_1, \phi_2: C\to A$ are two unital monomorphisms such that
$$
[\phi_1]=[\phi_2]\,\,\,{\rm in}\,\,\, KL(C, A)\andeqn
\\
\tau\circ \phi_1=\tau\circ \phi_2
\tforal \tau\in T(A).
$$
If also $K_1(A)=H_1(K_0(C), K_1(B))$ {\rm (}see {\rm \cite{Lnpams2})}, then there exists a sequence of unitaries $u_n\in U_0(A)$ such that
$$
\lim_{n\to\infty}{\rm ad}\, u_n\circ \phi_1(f)=\phi_2(f)\tforal f\in C.
$$

\end{thm}

\begin{proof}
This follows from 3.6 of \cite{Lnind} that the above statement holds without requiring $u_n$ in $U_0(A).$  The reason that
$u_n$ can be taken in $U_0(A)$ is given in \cite{Lnpams2}--see the proof
of 12.4 of \cite{Lnajm}.

\end{proof}

\begin{lem}\label{U0}
Let $A$ be a unital simple \CA\, with stable rank one and let $F$
be a  finite dimensional \SCA. Suppose that there are two
monomorphisms $\phi_1, \phi_2: F\to A$ such that $\phi_1$ is not
unital and
$$
(\phi_1)_{*0}=(\phi_2)_{*0}.
$$
Then there exists a unitary $u\in U_0(A)$ such that
$$
u^*\phi_1(a)u=\phi_2(a),\,\,\,\rforal a\in F.
$$
\end{lem}

\begin{proof}
Since $A$ has stable rank one, it is well known that there is a
unitary $v\in U(A)$ such that
$$
v^*\phi_1(a)v=\phi_2(a)\,\,\rforal a\in F.
$$
Since $\phi_1$ is not unital and $A$ is of stable rank one,
neither is $\phi_2.$ Let $e=1-\phi_2(1_F).$ Then $e\not=0.$

Since $A$ is simple, $eAe$ is stably isomorphic to $A.$
Furthermore, the map $w\mapsto (1-q_1)+w$ is an isomorphism from
$K_1(q_1Aq_1)$ onto $K_1(A).$ Therefore, since $A$ has stable rank
one, there exists $w\in U(eAe)$ such that $[(1-e)+w]=[v]$ in
$K_1(A).$ Define
$$
u=v(1-e+w)^*.
$$
Then $u\in U_0(A).$  Moreover,
$$
u^*\phi_1(a)u=\phi_2(a)\,\,\,\rforal a\in F.
$$

\end{proof}


\begin{df}\label{df1}
Let $f: S^1\to S^1$ be a degree $k$ map ($k> 1$), i.e., a continuous map with the winding number $k.$
Following 4.2 of \cite{EG}, denote by $T_{II,k}=D^2\cup_fS^1$ the connected finite CW complex obtained by
attaching a $2$-cell $D^2$ to $S^1$ via the map $f.$

Let $g: S^2\to S^2$ be a degree $k$ ($k>1$) map. Let $T_{III, k}=D^3\cup_g S^2$ the connected finite CW complex obtained by attaching a $3$-cell $D^3$ to $S^2$ via the map $g.$

Let $C=\oplus_{i=1}^rC_i,$ where $C_i=P_iM_{k_i}(C(X_i))P_i,$ where $P_i\in M_{r_i}(C(X_i))$ is a
projection and $X_i$ is a point, $X_i=S^1,$  $X_i=T_{II, m_i}$ or $T_{III, M_i}.$
\end{df}

\begin{lem}\label{Bot}
Let $C=\oplus_{i=1}^r C_i$ be a unital \CA, where $C_i$ is as described in \ref{df1},  Let $\ep>0,$ let ${\cal F}\subset C$ be a finite subset. There is $\dt>0$ and a finite subset ${\cal G}\subset  C$ satisfying the following. Suppose that $A$ is a unital separable simple \CA\,  with $TR(A)=0$ and suppose that $\phi: C\to A$ be a unital monomorphism and $u\in U(A)$ is a unitary such that
\beq\label{Bot-1}
\|[u, \, \phi(g)]\|<\dt,
\eneq
then there exists a unitary $v\in U(A)$  such that
\beq\label{Bot-2}
\|[v, \phi(f)]\|<\ep\tforal f\in {\cal F}\andeqn {\rm Bott}(uv, \phi)=0,\\
{\rm or}\,\,\, {\rm Bott}(vu, \phi)=0.
\eneq
\end{lem}
This is a combination of 6.7, 6.8, 6.9 and 6.10 of \cite{Lnajm}.

\begin{lem}\label{nieqi}
Let $A$ be a unital simple \CA\, of tracial rank zero, let $C$ be another \CA\, such that $A\subset C$ with $1_C=1_A,$  and let $B_1=\oplus_{i=1}^r C_i,$ where $C_i$ is as described in \ref{df1} with $1_C=1_A.$ Let $B_2$ be a
finite dimensional \CA\, with $p=1_{B_2}.$
Suppose that $px=xp$ for all $x\in B_1,$
$x\mapsto (1-p)x(1-p)$ is injective and $pC_lp\not=0,$ $l=1,2,...,r.$ Let $B_3=B_2\oplus (1-p)B_1(1-p)$ and let $\phi_0: B_3\to A$ be the embedding.

Suppose that $U\in U(C)$ with $U^*aU\in A,$
$\tau(U^*aU)=\tau(a)$
for all $a\in A$ and $\tau\in T(A)$ and
\beq\label {neiqi-0}
[{\rm ad}U\circ \phi_0]=[\phi_0]\,\,\,{\rm in}\,\,\,KL(A, A).
\eneq

Let $\ep>0,$ ${\cal F}\subset B_1$ and let $m\ge 1$ be a finite subset.
Then there exists
$\dt>0$ and a finite subset ${\cal G}\subset B_1$ satisfying the following:

Suppose that $V\in U(A)$ with $V^*U^*1_{C_j}U V=1_{C_j},$ $1\le j\le r,$ such that
\beq\label{nieqi-1-}
\|V^*U^*aUV-a\|<\dt\tforal a\in {\cal G}
\eneq
and suppose that $V=V_1V_2\cdots V_m$ for some $m\ge 2$ with
\beq\label{nieqi-1}
\|V_i-1\|<{\pi\over{(m-1)}},\,\,\,i=1,2,...,m.
\eneq

Then there exist unitaries $W_1, W_2,...,W_m\in U(A)$ such that
\beq\label{neiqi-2}
\|W_i-1\|<{2\pi\over{m-1}}+\|V_i-1\|,\,\,\,i=1,2,...,m-1,\\
\|(W_1W_2\cdots W_l)b(W_1W_2\cdots W_l)^*-
(V_1V_2\cdots V_l)b(V_1V_2\cdots V_l)^*\|<\ep
\eneq
for all $b\in {\cal F}$ and
\beq\label{neiqi-3}
(W_1W_2\cdots W_m)^*U^*aU(W_1W_2\cdots W_m)=a
\eneq
for all $a\in B_2.$

\end{lem}

\begin{proof}
Let $E(l)$ be the identity of $C_l,$ $l=1,2,...,r.$
By the assumption,
\beq\label{neiqi-4}
V^*U^*E(l)UV=E(l), \,\,\, l=1,2,...,r.
\eneq
Note that $pE(l)=E(l)p\not=0,$ $l=1,2,...,r.$
Therefore, by considering each summand individually, without loss of generality,
 we may assume that $r=1.$   So for the rest of the proof, $C_1=B_1.$
 Denote by $\phi: C_1\to A$ the embedding and denote by $\phi': C_1\to (1-p)A(1-p)$ and $\phi'': C_1\to B_2\subset pAp$
 the \hm s defined by $\phi'(c)=(1-p)\phi(c)(1-p)$ and
 by $\phi''(c)=pcp$
  all $c\in C_1,$ respectively.



Therefore one has  $pB_1=B_1p=M_{k(1)}.$
Then $B_2=M_{k_1}(F)$ for some finite dimensional \CA\, $F.$
By \ref{U0}, there exists a unitary $W'\in U_0(A)$ such that
\beq\label{neiqi-5}
(W')^*V^*U^*bUVW'=b\tforal b\in B_2.
\eneq
In particular, for all $c\in C_1=B_1,$
\beq\label{neiqi-6}
(W')^*V^*U^*pcpUVW'=pcp,
\eneq
Choose $\dt_1>0$ and a finite subset ${\cal G}_1\subset B_1$
such that the following holds:
\beq\label{neiqi-n1}
{\rm Bott}(Z, \psi_1)={\rm Bott}(Z, \psi_2)
\eneq
for any unitary $Z\in U(A)$ and any unital \hm s $\psi_1, \psi_2: B_1\to A,$ whenever
$$
\|\psi_1(f)-\psi_2(f)\|<\dt_1\tforal f\in {\cal G}_1,
$$
provided that both sides in (\ref{neiqi-n1}) are defined
(see \cite{Lnmemoir}).

Our $\dt$ and ${\cal G}$ will be chosen later, but, at least,
we will choose $\dt<\dt_1$ and ${\cal G}\supset {\cal G}_1.$

For any given $\dt_2>0$  and  finite subset ${\cal G}_2\subset B_1,$
let  $\dt_3>0$ be as in \ref{Bot} (in place of $\dt$) and let
${\cal G}_3$ be a finite subset as in \ref{Bot} (in place of ${\cal G}_2$) associated with
$\dt/4$ (in place of $\ep$) and ${\cal G}$ (in place of ${\cal F}$).
By 3.4 of \cite{Lncd} (see also 3.6 of \cite{Lnind}), from the assumption (\ref{neiqi-0}) and
the assumption that $\tau(U^*xU)=\tau(x)$ for all $x\in A$ and $\tau\in T(A),$ there is a unitary $U'\in U(A)$ such that
\beq\label{neiqi-n2}
&&\|(U')^*cU'-U^*cU\|<\dt_3/2<\dt_1/4\tforal c\in {\cal G}_3
\andeqn\\\label{neiqi-n2+}
&&(U')^*bU'=U^*bU\tforal b\in B_2.
\eneq
Moreover,  by applying \ref{Bot},
if $\|V^*U^*cUV-c\|<\dt_3/2$ for all $c\in {\cal G}_3,$
we can obtain another unitary $U''\in U(A)$ such that $(U'')^*bU''=b$ for all $b\in B_2,$ 
\beq\label{neiqi-nn1}
\|U''c-cU''\|<\dt_2\tforal c\in {\cal G}_2\andeqn {\rm Bott}(VU'U'',\phi)=0.
\eneq

By choosing $\dt<\dt_3,$ ${\cal G}\supset {\cal G}_3$ and  replacing $U'$ by $U''U',$ simplifying the notation,  we may simply assume (omitting $\dt_2$ and ${\cal G}_2$) that
\beq\label{neiqi-nn1}
{\rm Bott}(VU', \phi)=0,\,\,\,\|(U')^*cU'-U^*cU\|<\dt/4\tforal c\in {\cal G}\\\label{neiqi-nn2}
\andeqn (U')^*bU'=U^*bU\tforal b\in B_2
\eneq
By (\ref{neiqi-0}) and (\ref{neiqi-6}),
\beq\label{neiqi-7}
[{\rm ad} UVW'\circ \phi']=[\phi'].
\eneq
Then, by applying \cite{Lncd} (see also 3.6 of \cite{Lnind}), if $\dt>0$ (we may assume that $\dt<\ep/4m$) and a finite subset ${\cal G}\subset B_1$ are given,
there is a unitary $W''\in (1-p)A(1-p)$ such that
\beq\label{neiqi-8}
\|(W'')^*(W')^*V^*U^*(1-p)c(1-p)UVW'W''-(1-p)c(1-p)\|<\dt/8\tforal c\in {\cal G}.
\eneq
It follows that
\beq\label{neiqi-n3}
\|(W'')^*(W')^*V^*(U')^*(1-p)c(1-p)U'VW'W''-(1-p)c(1-p)\|<3\dt/8
\eneq
for all $c\in {\cal G}.$

For the monomorphism $\phi,$ we choose $\eta>0$ and a finite subset ${\cal G}_2\in B_1$
so that 17.5 (see also 8.4)  of \cite{Lnmemoir} can be applied
for $\ep/2m$ (in place of $\ep$) and ${\cal F}.$  We may assume that ${\cal F}\subset {\cal G}_2.$

Then, by choosing large ${\cal G}$ and small
$\dt,$  by \ref{Bot}, there is a unitary $W_0\in U((1-p)A(1-p))$
such that
\beq\label{neiqi-9}
&&\|W_0(1-p)c(1-p)-(1-p)c(1-p)W_0\|<\eta/4\tforal {\cal G}_1\andeqn\\\label{neiqi-9+1}
&&{\rm Bott}(((1-p)U'VW'W''W_0), \phi'))=0.
\eneq
Moreover, we may assume that $\dt<\eta/4$ and ${\cal G}\supset {\cal G}_2.$
Keep in mind that we also assume that (\ref{nieqi-1-}) holds for the above $\dt$ and ${\cal G}.$
Put $W=W'(W''\oplus p)(W_0\oplus p).$ Then,
\beq\label{neiqi-9+}
\hspace{-0.8in}W^*V^*U^*cUV^*W &=& W^*V^*U^*(1-p)c(1-p)UVW+W^*V^*U^*qcqUVW\\
&\approx_{\dt+\eta/4}& (1-p)c(1-p)+pcp=c\approx_{\dt} V^*U^*cUV\\
&& \tforal c\in {\cal G}_1\andeqn\\\label{neiqi-n9}
W^*V^*U^*bUVW&=& p(W')^*V^*U^*bUVW'p=pbp=b\tforal b\in B_2.
\eneq
In particular, (by applying (\ref{neiqi-nn2})),
\beq\label{neiqi-n10}
pU'VW=U'VWp\andeqn (pU'VWp)b=bpU'VWp\tforal b\in B_2.
\eneq
Let
\beq\label{neiqi-n12}
W_{00}=(pU'VWp)^*\in U(pAp).
\eneq
By replacing $W$ by $W((1-p)\oplus W_{00}),$ we may assume that $pU'VWp=p.$
Then, by (\ref{neiqi-9+1}), we compute that
\beq\label{neiqi-n5}
{\rm Bott}(U'VW, \phi)&=&{\rm Bott}((1-p)U'VW, \phi')+{\rm Bott}(pU'VW, \phi'')=0.
\eneq
But,  by (\ref{neiqi-nn1}), we also have that
\beq\label{neiqi-n6}
0={\rm Bott}(U'VW, \phi)&=&{\rm Bott}(U'V, \phi)+{\rm Bott}(W, \phi)\\
={\rm Bott}(W, \phi).
\eneq
It follows from  the choice of $\dt_1$ and ${\cal G}_1$  and
(\ref{neiqi-n1}) that
\beq\label{neiqi-n7}
{\rm Bott}(W, {\rm ad}\, U'V\circ \phi)=0.
\eneq
It follows from 17.5 (see also 8.4) of \cite{Lnmemoir} that there are $Z_1,Z_2,..., Z_m\in U(A)$ such that
\beq\label{neiqi-10}
&&\hspace{-0.4in}W=Z_1Z_2\cdots Z_m,\,\,\,\|W_i(V^*U^*bUV)-(V^*U^*bUV)W_i\|<\ep/2m\tforal b\in  {\cal F}\\
&&\andeqn \|Z_i-1\|<2\pi/(m-1),\,\,\,i=1,2,...,m.
\eneq
So, by (\ref{nieqi-1-}), (note also we assume that $\dt<\ep/4m$),
\beq\label{neiqi-10+}
\|W_ib-bW_i\|<\ep/m, i=1,2,...,m.
\eneq
Now define
\beq\label{neiqi-11}
W_1&=&V_1Z_1,\, \, X_2=Z_1^*V_2Z_1, Z_2,\,...\\
W_j&=&(Z_1Z_2\cdots Z_{j-1})^*V_j(Z_1Z_2\cdots Z_{j-1}),..., \\
W_m&=&(Z_1Z_2\cdots Z_{m-1})^*V_m(Z_1Z_2\cdots Z_{m-1}).
\eneq
We estimate that
\beq\label{neiqi-12}
\|W_j-1\| &\le & \|(Z_1Z_2\cdots Z_{j-1})^*(V_j-1)(Z_1Z_2\cdots Z_{j-1})Z_l\|+\|Z_j-1\|\\
&<& \|V_{j-1}-1\|+{2\pi\over{m-1}},\,\,\,l=1,2,...
\eneq
Note that
\beq\label{neiqi-13}
W_1W_2\cdots X_j=V_1V_2\cdots V_jZ_1Z_2\cdots Z_j,\,\,\,j=1,2,...,m
\eneq
For each $a\in B_2,$
\beq\label{neiqi-14}
(W_1W_2\cdots W_m)^*U^*bU(W_1W_2\cdots W_m) &=&
W^*V^*U^*bUVW=b\tforal b\in B_2.
\eneq
Moreover, for $c\in {\cal F},$
\beq\label{neiqi-15}
&&\hspace{-0.8in}(W_1W_2\cdots W_j)c(W_1W_2\cdots W_j)^* \\
&=&V_1V_2\cdots V_jZ_1Z_2\cdots Z_jc(Z_1Z_2\cdots Z_j)^*(V_1V_2\cdots V_j)^*\\
&\approx_{j\ep/m}& (V_1V_2\cdots V_j)c(V_1V_2\cdots V_j)^*,\,\,\,j=1,2,...,m.
\eneq

\end{proof}

\begin{lem}\label{Ad4}
Let $B$ be a finite dimensional \SCA\, of a unital simple \CA\,
$A.$ Then for any  $\dt>0$ there exists $\sigma>0$ and a finite
subset ${\cal G}\subset A$ satisfying the following: If
$$
\|pf-fp\|<\sigma
$$
for all $f\in {\cal G},$ then there is an monomorphism $\phi:B\to
pAp$ or $\phi: B\to 1_BA1_B$ such that
$$
\|pbp-\phi(b)\|<\dt\|b\|\rforal b\in B.
$$
\end{lem}

\begin{proof}
Write $B=M_{r(1)}\oplus M_{r(2)}\oplus\cdots\oplus M_{r(l)}.$ Let
$e_{i,j}^{(s)}\in B$ be a system of matrix units for
$M_{r(s)},$ $s=1,2,...,l.$ Since $A$ is simple, it is easy to
obtain, for each $i,$ elements
$v_{s,i,1},v_{s,i,2},...,v_{s,i,m(s)}$ in $B$ with
$\|v_{s,i,k}\|\le 1$ such that
$$
\sum_{k=1}^{m(s)}v_{s,i,k}^*e_{i,i}^{(s)}v_{s,i,k}=1_A,
s=1,2,...,l.
$$
Let ${\cal F}_0=\{e_{i,j}^{(s)}:1\le s\le l\}\cup\{v_{s,i,k}: 1\le
i\le R(s), 1\le k\le m(s), 1\le s\le l, \}.$ Let $J=\max\{m(s):
s=1,2,...,l\}.$ Let $\eta>0$ to be determined. Suppose that
$$
\|pb-bp\|<(1/2J)^2\min\{1/4,\dt/4\}\rforal b\in {\cal F}_0.
$$
Then
$$
\|\sum_{k=1}^{m(s)}pv_{s,i,k}^*pe_{i,i}^{(s)}pv_{s,i,k}p-p\|<(1/2J)\min\{1/4,\eta/4\},
s=1,2,...,l.
$$
It follows that
$$
\|pe_{i,i}^{(s)}p\|>1/2J, i=1,2,...,R(s), s=1,2,...,l.
$$
Put $a=pe_{i,i}^{(s)}p.$ We claim that, in fact, $\|a\|\ge 1/2.$
Otherwise $\|a\|<1/2.$ We have
\begin{eqnarray}\nonumber
\|a-a^2\|=\|pe_{i,i}^{(s)}p-(pe_{i,i}^{(s)}pe_{i,i}^{(s)}p)\|<(1/2J)^2\min\{1/4,\eta/4\}.
\end{eqnarray}
Applying the spectral theorem, if $t=\|a\|,$ we have  that
$$
|t-t^2|<(1/2J)^2\min\{1/4,\eta/4\}.
$$
It follows that (since $|t|=\|a\|>1/2J$)
$$
|1-t|<(1/2J)(1/4)\le 1/8.
$$
It is impossible unless $t\ge 1/2.$

It then follows from 2.5.5 of \cite{Lnbk} that there is a nonzero
projection $q_{i,i}^{(s)}\in pAp$ such that
$$
\|pe_{i,i}^{(s)}p-q_{1,1}^{(s)}\|<(1/2J)^2\min\{1/2,\eta/2\}.
$$
The rest of proof is standard and follows from the argument in
section 2.5 of \cite{Lnbk} and  2.3 of \cite{Lntqe}.
\end{proof}

\begin{lem}\label{Ad5}
Let $A$ be a unital \CA\, and let $V, V_1,V_2,...,V_m\in U_0(A)$
be unitaries such that $V=V_1V_2\cdots V_m.$

Then, for any $\dt>0$ and any nonzero projection $p\in A$ there is
$\eta=\eta(\dt,m)>0$ (which does not depend on $A$) satisfying the following:
if
$$
\|pV-Vp\|<\eta\andeqn \|pV_i-V_ip\|<\eta,\,\,\,i=1,2,...,m,
$$
then there exist unitaries $W, W_i\in pAp$ such that
$$
\|W-pVp\|< \dt,\,\,\, \|W_i-pV_ip\|<\dt \andeqn W=W_1W_2\cdots
W_m.
$$

\end{lem}

The proof of the above is standard (see for example section 2.5 of
\cite{Lnbk}).

\begin{lem}\label{UD}
Let $A$ be a unital  \CA.
Let $k>1$ be an integer. Suppose that there are $k$
mutually orthogonal projections $e_1,e_2,...,e_k$ in $A$ and a
unitary $u\in U(A)$ such that
$$
u^*e_iu=e_{i+1},\,\,\,i=1,2,...,k\,\,\,\andeqn e_{k+1}=e_1.
$$
Suppose that $B$ is a finite dimensional \SCA\, in $e_1Ae_1$  and
$z\in U_0(e_1Ae_1)$ such that $z^*(u^k)^*bu^kz=b$ for all $b\in
B.$ Suppose that $z=z_1z_2\cdots z_{k-1},$ where $z_i\in
U(e_1Ae_1),$ $i=1,2,...,k-1$ nd $z_k=1.$ Define
$$
w=\sum_{i=1}^ke_iu^{k+1-i}z_i(u^{k-i})^*+(1-\sum_{i=1}^ke_i)u.
$$
Then
\begin{align}\nonumber
&\|w-u\|\le \max\{\|z_i-1\|:1\le i\le k-1\}\\
&(w^i)^* e_1(w^i)=
e_{i+1},\,\,\,\,\,\,i=1,2,...,k-1\andeqn\\\nonumber
&(w^{k})^*bw^{k}=b\,\,\,\,\,\,\,\,\,\rforal b\in B.\nonumber
\end{align}
\end{lem}

\begin{proof}
It is easy to verify that $w$ is a unitary. One estimates that
$$
\|w-u\|\le \max\{\|z_i-1\|:1\le i\le k-1\}.
$$
Note that
\begin{eqnarray}\nonumber
(w^i)^*e_1w^i&=&u^{k-i}(z_1z_2\cdots
z_i)^*(u^k)^*e_1u^k(z_1z_2\cdots z_i)(u^{k-i})^*\\\nonumber
&=&u^{k-i}e_1(u^{k-i})^*=e_{i+1},\,\,\,\,\,\,\,\,\,\hspace{0.3in}
i=1,2,...,k.
\end{eqnarray}
Moreover,
\begin{eqnarray}\nonumber
e_1w^k&=&e_1u^k(z_1z_2\cdots z_{k-1})\\
&=&e_1u^kz.
\end{eqnarray}
Then
$$
(w^k)^*bw^k=z^*(u^k)^*bu^kz=b\rforal b\in B.
$$
 The lemma then follows.

\end{proof}

\begin{lem}\label{GR}
Let $A$ be a unital \CA. Suppose that $e_1,e_2,...,e_n$ are $n$
mutually orthogonal projections and $u\in U(A)$ is a unitary such
that
$$
u^*e_iu=e_{i+1},\,\,\, i=1,2,...,n\andeqn e_{n+1}=e_1.
$$
Then $\{e_i: i=1,2,...,n\},$ $\{ e_iue_{i+1}: i=1,2,...,n-1\},$
$\{pup\},$ where $p=\sum_{i=1}^ne_i,$ generate a \SCA\, which is
isomorphic to $C(X)\otimes M_n,$ where $X$ is a compact subset of
$S^1.$
\end{lem}

\begin{proof}
It is standard to check that $\{e_i:i=1,2,...,n\}$ and
$\{e_iue_{i+1}: i=1,2,...,n-1\}$ generate a \SCA\, which is
isomorphic to $M_n.$ Let $z=e_1u^{n}e_1.$ Since $e_1u^n=u^ne_1,$
$z$ is a unitary in $e_1Ae_1.$ Suppose that $sp(z)=X\subset S^1.$
Then $\{e_i:i=1,2,...,n\},$ $\{e_iue_{i+1}: i=1,2,...,n-1\}$ and
$z$ generate a \SCA\, $B$ which is isomorphic to $C(X)\otimes
M_n.$ Note
$$
e_1u^{n-1}e_n=e_1ue_2ue_3\cdots e_{n-1}ue_n\in B.
$$
We also have that
$$
e_nue_1=(e_1u^{n-1}e_n)^*z\in B.
$$
So
$$
pup=\sum_{i=1}^{n-1} e_iue_{i+1}+ e_nue_1\in B.
$$

On the other hand, the \SCA\, generated by $\{e_i:i=1,2,...,n\},$
$\{e_iue_{i+1}: i=1,2,...,n-1\}$ and $\{pup\}$ contains $e_nue_1$
as well as $z.$ This proves the lemma.
\end{proof}

\section{Tracial rank zero}

\begin{lem}\label{Ad6}
Let $A$ be a unital separable simple \CA\, with tracial rank zero
and let $G\subset K_0(A)$ be a subgroup such that $\rho_A(G)$ is
dense in $\rho_A(K_0(A)).$ Then, for any $\ep>0$ any finite subset
${\cal F}\subset A$ and any nonzero positive element $b\in A_+,$
there exists a projection $p\in A$ and a finite dimensional \SCA\,
$B\subset A$ with $1_B=p$ and $[e]\in G$ for all projections $e\in
B$ such that

{\rm (1)} $\|pa-ap\|<\ep\rforal a\in {\cal F},$

{\rm (2)} $pap\in_{\ep} B$ for all $a\in {\cal F}$ and

{\rm (3)} $[1-p]\le [b].$
\end{lem}

\begin{proof}
Fix any  $\ep>0$ and any finite subset ${\cal F}\subset A$ and any
nonzero positive element $b\in A_+.$ There are mutually orthogonal
nonzero projections $r_1, r_2\in {\overline{bAb}}.$ Since
$TR(A)=0,$ there is a projection $q\in A$ and there is a finite
dimensional \SCA\, $B_1\subset A$ with $1_{B_1}=q$ such that

{\rm (1)} $\|qa-aq\|<\ep\rforal a\in {\cal F},$

{\rm (2)} $qaq\in_{\ep} B_1$ for all $a\in {\cal F}$ and

{\rm (3)} $[1-q]\le [r_1].$

Write $B_1=M_{R(1)}\oplus M_{R(2)}\oplus \cdots \oplus M_{R(k)}.$
Let $\{e_{i,j}^{(l)}\}$ be a system of matrix units,
$l=1,2,...,k.$ Put $\dt=\inf\{\tau(r_2): \tau\in T(A)\}.$ Note
that $\dt>0.$ Since $TR(A)=0$ and $\rho_A(G)$ is dense in
$\rho_A(K_0(A)),$ there is, for each $l,$ a nonzero projection
$d_l\le  e_{1,1}^{(l)}$ such that $[d_l]\in G$ and
$$
\tau(e_{1,1}^{(l)}-d_l)<{\dt\cdot\inf\{\tau(e_{1,1}^{(l)}:\tau\in
T(A)\}\over{R(l)}}\rforal \tau\in T(A),\,\,\,l=1,2,...,k.
$$
Define
$$
p_l=\sum_{i=1}^{R(l)}e_{i,1}^{(l)}d_le_{1,l}^{(l)},\,\,\,l=1,2,...,k.
$$
In other words, $p_l$ has the form ${\rm diag}(\overbrace{d_l,
d_l,\cdots, d_l}^{R(l)}).$ It follows that $[p_l]\in G$ and $p_l$
commutes with every element in $B_1.$ We compute that
$$
\tau(\sum_{i=1}^{R(l)}e_{i,i}^{(l)}-p_l)<\dt\cdot\inf\{\tau(e_{1,1}^{(l)}:\tau\in
T(A)\}\rforal \tau\in T(A),\,\,\, l=1,2,...,k.
$$

Define
$$
p=\sum_{l=1}^k p_l\andeqn B=pB_1p.
$$
It follows that
$$
\tau(q-p)<\dt<\tau(r_2)\rforal \tau\in T(A).
$$
Therefore, by \cite{Lntr0}, $[p-q]\le [r_2].$ On the other hand,
since $[d_l]\in G,$ we see that $[e]\in G$ for all projections
$e\in B.$

Moreover, we have the following:

(i) $\|pa-ap\|<\ep\rforal a\in {\cal F},$

(ii) $pap\in_{\ep} B$ for all $a\in {\cal F}$ and

(iii) $[1-p]\le [r_1]+[r_2]\le [b].$

\end{proof}

\begin{lem}\label{one+0}
Let $A$ be a unital separable simple \CA\, with real rank zero.
Suppose that $A$ has the following property:

For any $\ep>0,$ any finite subset ${\cal F}\subset A$ and any
nonzero positive element $b\in A_+,$ there exists a projection
$p\in A$ and a \SCA\, $D\subset A$ with $1_D=p$ and $D\cong
\oplus_{j=1}^N C(X_j)\otimes F_j,$ where each $F_j$ is a finite
dimensional \SCA\, and $X_j\subset S^1$ is a compact subset, such
that

{\rm (1)} $\|pa-ap\|<\ep\rforal a\in {\cal F},$

{\rm (2)} $pap\in_{\ep} D$ for all $a\in {\cal F}$ and

{\rm (3)} $[1-p]\le [b].$

Then $A$ has tracial rank zero.
\end{lem}

\begin{proof}
This is known. We sketch the proof as follows: Since $A$ has real
rank zero, $C(X_j)\otimes F_j$ is approximated pointwise in norm
by a finite dimensional \SCA, if $X_j\not=S^1.$ It is then easy to
see , with (1), (2) and (3) above, that $TR(A)\le 1$ (by replacing
$\ep$ by $2\ep$ and replacing $C(X_j)\otimes F_j$ by a finite
dimensional \SCA, if $X_j\not=S^1$) (see also 6.13 of \cite{Lntr0}
).

It follows from (b) of Theorem 7.1 of \cite{Lntr0} that $A$ is TAI
(see the definition of TAI, for example, right above 7.1 of
\cite{Lntr0}). Since $A$ has real rank zero, any \SCA\, with the
form $M_m(C(I)),$ where $I=[0,1],$ is approximated by finite
dimensional \SCA s. It follows that $TR(A)=0.$
\end{proof}

\begin{rem}

 An alternative proof (without using Theorem 7.1 of \cite{Lntr0}) is described below:

First show that $A$ has stable rank one (the same proof as that
$TR(A)\le 1$ implies that $A$ has stable rank one). Or, as above,
first show that $TR(A)\le 1.$ Then one  concludes that $A$ has
stable rank one.

Let  $u\in A$ be a unitary such that $sp(u)=S^1.$ So we have a
monomorphism $h_1: C(S^1)\to A.$ Let ${\cal G}\subset C(S^1)$ be
any finite subset and $\ep>0.$ Let $\dt>0$ be also given. For any
$\eta>0$ and  for any projection $r\in A,$ there is a projection
$e$ and a unitary $v\in (1-e)A(1-e)$ such that $[e]\le [r]$ and
$$
\|u-(e+v)\|<\eta.
$$
Write $e=e_1+e_2,$ where $e_1, e_2$ are two non-zero mutually
orthogonal projections. Since $A$ has stable rank one, choose
$v_1\in e_1Ae_1$ and $v_2\in e_2Ae_2$ such that
$[v_2]=[v^*]=[u^*]$ and $[v_1]=[v]=[u]$ in $K_1(A).$ Put
$w=v_1+v_2+v.$ Verify that $[u]=[w]$ in $K_1(A)$ and verify that
$$
|\tau(f(u))-\tau(f(w))|<\dt \rforal \tau\in T(A) \andeqn \rforal f\in {\cal G}
$$
if $\eta$ is sufficiently small and $\sup\{\tau(r): \tau\in
T(A)\}$ is sufficiently small. It follows from, for example,
Theorem 3.3 of \cite{Lncd}, that there is a unitary $z\in A$ such
that
$$
\|u-z^*wz\|<\ep
$$
if $\dt$ is sufficiently small.

In other words, there is a unitary $U\in A$ with the form
$U=U_1+U_2,$ where $U_1\in qAq$ and $U_2\in (1-q)A(1-q)$ are
unitaries for which  $[q]\le [e_1]\le [r]$ and $[U_2]=0$ in
$K_1(A),$ such that
$$
\|u-U\|<\ep.
$$

Since $A$ is assumed to have real rank zero, $(1-q)A(1-q)$ has
real rank zero. It follows from \cite{Lnfu} that $U_2$ can be
approximated in norm by unitaries with finite spectrum. From this
and together with (1),(2) and (3) above, one easily sees that
$TR(A)=0.$

\end{rem}

\begin{thm}\label{TM}
Let $A$ be a unital separable amenable simple \CA\, with $TR(A)=0$
which satisfies the UCT. Suppose that
$\alpha\in Aut(A)$ has the tracial cyclic Rokhlin property.
Suppose also that  there is an integer $J\ge 1$ such that
$[\alpha^J]=[{\rm id}_A]$ in $KL(A,A).$
 Then
$TR(A\rtimes_{\alpha}\Z)=0.$

\end{thm}

\begin{proof}
We first note that, by \cite{K1}, $A\rtimes_{\alpha}{\mathbb Z}$
is a unital simple \CA. By \cite{OP}, $A\rtimes_\af\Z$ has real rank zero.  We will show that $(A\rtimes_\af\Z)\otimes Q$ has tracial rank one, where $Q$ is the UHF-algebra
 with $(K_0(Q), [1_Q])=(\Q, 1).$ Let $u\in A\rtimes_\af\Z$ be a unitary
which implements
$\alpha,$ i.e., $\alpha(a)=u^*au$ for all $a\in A.$
Put $B=A\otimes Q.$ 
Note that $(A\rtimes_\af\Z)\otimes Q$ is generated by $A\otimes Q$ and elements of the form 
$u\otimes a$ for $a\in Q.$  One can identify $1_A$ and $1_{A\rtimes \af\Z}.$ Since $A\otimes Q$ contains elements $1_A\otimes a$ ($a\in Q$), $(A\rtimes_\af\Z)\otimes Q$ is also generated by $A\otimes Q$ and $u\otimes 1_Q.$ 
By identifying $u$ with $u\otimes 1_Q,$ 
$(A\rtimes_\af\Z)\otimes Q$ is generated by $B$ and $u.$ So, in what follows, we will identify $u$ with $u\otimes 1_Q.$ 
We will first show that $TR((A\rtimes_\af \Z)\otimes Q)=0.$

To this end, let $1>\ep > 0$ and ${\cal F} \subset (A \rtimes_\alpha{\mathbb Z})\otimes Q$
be a finite set. To simplify notation, without loss of generality,
we may assume that
$$
{\cal F} = {\cal F}_0\cup \{u\},
$$
where ${\cal F}_0\subset  B$ is a finite subset of the unit ball
which contains $1_{B}.$

 Choose an integer $k$ which is a
multiple of $J$ such that $8\pi/(k-2)<\ep/256.$
Put
\begin{eqnarray}\label{eI}
{\cal F}_1={\cal F}_0\cup \{u^ia(u^*)^i: a\in {\cal F}_0, -k\le
i\le k \}.
\end{eqnarray}

Fix $b_0\in ((A \rtimes_\alpha{\mathbb Z})\otimes Q)_+\backslash \{0\}.$ It
follows from Theorem 4 of \cite{JO}  that $A\rtimes_{\alpha}\Z$
has property (SP). Thus there is a nonzero projection $r_{00}$ in
the hereditary \SCA\, of $(A\rtimes_{\alpha}\Z)\otimes Q$ generated by $b_0.$

Let $r_1',r_2'\in B$ be  nonzero mutually orthogonal
projections. Since $A\rtimes_{\alpha}\Z$ is simple, it follows
from \cite{C} (see also (2) of 3.5.6 and 3.5.7 in \cite{Lnbk})
that there are nonzero projections $r_i\le r_i',$ such that
$[r_i]\le [r_{00}],$ $i=1,2.$

 Let $\eta_0={\ep\over{256k^4}}.$  Since $TR(A)=0$ and satisfies the UCT,  by the classification theorem of \cite{Lnduke} and \cite{EG},
 without loss of generality, we may assume that
 ${\cal F}_1\subset B_1\subset B,$ where
 $B_1=\oplus_{i=1}^rC_i$ and $C_i$ has the form described in
 \ref{df1}. Let $\dt>0$ and ${\cal G}\subset B_1$ be a finite subset in \ref{nieqi} corresponding to
 $\eta_0$ (in place of $\ep$) and ${\cal F}_1$ (in place of ${\cal F}$).
 Since $TR(A)=0$ and $[\af^k]=[{\rm ad}_A]$ in $KL(A,A),$
 \beq\label{TMn-1}
 \tau((u^k)^*au^k)=\tau(a)\tforal a\in B\andeqn \tforal \tau\in T(B).
 \eneq
 (Recall that we identify $u$ with $u\otimes 1_Q$).
 By the assumption and by applying 3.4 of \cite{Lncd}
 (see also 3.6 of \cite{Lnind}), there is a unitary $w\in U(B)$ such that
 \beq\label{TMn-2}
 w^*(u^k)^*1_{C_j}u^kw=1_{C_j},\,\,\,j=1,2,...,r\andeqn\\
 w^*(u^k)^*cu^kw\approx_{\dt/2} c\rforal c\in {\cal G}.
 \eneq
 Since $K_1(B)$ is divisible, $H_1(K_0(B_1), K_1(B))=K_1(B).$ From \ref{CDT}, we may assume that $w\in U_0(B).$ There are unitaries $w_1,w_2,...,w_{k-1}\in U(B)$ such that
 \beq\label{TMn-3}
 \|w_i-1\|<\pi/(k-1)\andeqn w=w_1w_2\cdots w_k.
 \eneq
Define
$$
{\cal F}_2=\{u^{-i}bu^i: b\in {\cal F}_1\cup {\cal G}, -k\le i\le k\}\andeqn
$$
\beq
&&\hspace{-0.9in}{\cal F}_3=\{u^i(w_{i_1}w_{i_1+1}\cdots
w_{i_1+l})^{i_2}a((w_{i_1}w_{i_1+1}\cdots w_{i_1+l})^{i_2})^* u^i:\\
&&\hspace{0.3in}b\in {\cal F}_2\cup {\cal G}, -k\le i,i_1,i_2\le k\}\cup\{w, w_j: 1\le j\le k\}.
\eneq
Let $\eta_1=\min\{\eta_0/4k, \dt/4k\}.$

We will again apply the classification theorem (\cite{Lnduke} and \cite{EG}).  In particular,
$A$ is an inductive limit of \CA s in \ref{df1} with injective connecting maps and with real rank zero. We may assume that there
is a projection $p_1\in B$ such that $p_1x=xp_1$ for all $x\in B_1$ and a finite dimensional
\SCA\, $B_2$ with $1_{B_2}=p_1$ such that $p_1x\in B_2$ for all $x\in B_1$ and $x\mapsto (1-p_1)x(1-p_1)$ is injective on $B_1,$ and

(a) $\|p_1f-fp_1\|<{\eta_1\over{64k}}$ for all $f\in {\cal F}_3,$

(b) $p_1fp_1\in_{{\eta_1\over{64k}}} B_2$  for all $f\in {\cal
F}_3$ and

(c) $[1-p_1]\le [r_1].$

Let $\phi': B_1\to (1-p_1)A(1-p_1).$
Put $B_3=B_2\oplus \phi'(B_1).$
It follows from
\ref{nieqi} that there exist unitaries $W, W_1,W_2,...,W_{k-1}\in B$
such that
\begin{eqnarray}\label{eMps1}
&&W^*(u^k)^*au^kW=a\rforal a\in B_2,\\\label{eMps2}
&&W=W_1W_2\cdots
W_{k-1},\,\,\,\|W_i-1\|<2\pi/(k-2)+\pi/(k-2)\andeqn
\\\label{eMps3}
&&(W_1W_2\cdots W_l)b(W_1W_2\cdots W_l)^*\approx_{\eta_0}(w_1w_2\cdots
w_l)b(w_1w_2\cdots w_l)^*
\end{eqnarray}
for all $b\in {\cal F}_1.$  Set
$$
{\cal F}_4={\cal F}_3\cup\{u^i(W_1W_2\cdots W_{i_1})^*u^ja
(u^i(W_1W_2\cdots W_{i_1})^*u^j)^*:a\in {\cal F}_4: -k\le i,i_1, j\le k\}.
$$
Let $\sigma_1>0$ and ${\cal G}_1$ be associated with $B_2$ and
$\eta_1/2$ in \ref{Ad4}. Let $\eta_2=\eta(\eta_1/2, k)$ be as
 in \ref{Ad5}. Let $\eta_3=\min\{\eta_2,\sigma_1, \eta_0\}.$ Define ${\cal
F}_5={\cal F}_4\cup {\cal G}_1.$

Since $\alpha$ has the tracial cyclic Rokhlin property,
 there exist
projections $e_1, e_2,...,e_k\in B$ such that
\begin{itemize}
\item[(1)] $\|\alpha(e_i) - e_{i+1}\| < {\eta_3\over{64k}}$ for $1 \leq i \leq
k$ ($e_{k+1} = e_1$)
\item[(2)] $\|e_ia-ae_i\| < {\eta_3\over{64k}}$ for $a\in
{\cal F}_5$
 \item[(3)] $[1 - \sum_{i=1}^ke_i] \leq [r_2]$
\end{itemize}
($e_i$ has the form $e_i\otimes 1_Q$).

Set $p = \sum_{i=1}^ke_i$. From (1) above, one estimates that
$$
\begin{array}{ll}
\|up - pu\| &= \|\sum_{i=1}^kue_{i+1} - \sum_{i=1}^ke_{i}u\|\\
&= \sum_{i=1}^k\|ue_{i+1} - e_{i}u\| =\sum_{i=1}^k\|ue_{i+1} -uu^*
e_{i}u\| < {\eta_3\over{64}}
\end{array}
$$

By (1) above, one sees that there is a unitary $v\in B$ such that
\begin{eqnarray}\label{eII}
\|v-1\|<{\eta_3\over{32k}}\andeqn
vu^*e_iuv^*=e_{i+1},\,\,\,i=1,2,...,k.
\end{eqnarray}
 Set
$u_1=v^*u.$ Then
$$
u_1^*e_iu_1=e_{i+1},\,\,\,i=1,2,...k\andeqn e_{k+1}=e_1.
$$
In particular,
$$
u_1^ke_1=e_1u_1^k.
$$

For any $a\in {\cal F}_5\cap B_2$ (since $W\in {\cal F}_5$),
$$
e_1W^*e_1(u_1^k)^*e_1ae_1u_1^ke_1We_1\approx_{\eta_3/16k} e_1ae_1.
$$

By \ref{Ad4}, there is a monomorphism $\phi_1: B_2\to e_1Be_1$
such that
\begin{eqnarray}\label{eMp1}
\|\phi_1(a)-e_1ae_1\|<{\eta_1\over{2}}\|a\|\rforal a\in B_2.
\end{eqnarray}

By applying \ref{Ad5}, and using (\ref{eMp1}), we obtain unitaries
$x, x_1, x_2,...,x_{k-1}\in U_0( e_1Be_1)$ such that
\begin{eqnarray}\label{eMp1+}
\|x-e_1We_1\|<\eta_1/2,\,\,\, \|x_i-e_1W_ie_1\|<\eta_1/2
\end{eqnarray}
$$
x=x_1x_2\cdots x_{k-1} \andeqn x^*(u_1)^*au_1x=a
$$
for all $a\in \phi_1(B_2).$









Let $Z=\sum_{i=1}^ke_iu_1^{k+1-i}x_i(u_1^{k-i})^*+(1-p)u_1.$ Define
$B_4=\phi_1(B_2).$ As in \ref{UD}, by (\ref{eMps2}),
\begin{eqnarray}\label{eMps4}
\|Z-u_1\|<\eta_1/2+3\pi/(k-2), \,\,(Z^k)^*bZ^k=b\rforal b\in
B_4
\end{eqnarray}
and $(Z^i)^*e_1Z^i= e_{i+1},$ $i=1,2,...,k$ ($e_{k+1}=e_1$).

Write $B_4=F_1\oplus F_2\oplus \cdots\oplus F_N$ and let $\{c_{is}^{(j)}\}$
be the matrix units for $F_j,$ $j=1,2,...,N,$ where $F_j=M_{R(j)},$
and put $q=1_{B_4}.$

Define $D_0=B_4\oplus \oplus_{i=1}^{k-1}Z^{i*}B_4Z^i$ and $D_1$
the \SCA\, generated by $B_4$ and $c_{ss}^{(j)}Z^i,$
$s=1,2,...,R(j),$ $j=1,2,...,N$ and $i=0,1,2,...,k-1.$ Then
$D_1\cong B_4\otimes M_k$ and $D_1\supset D_0.$

Define $q_{ss}^{(j)}=\sum_{i=0}^{k-1}Z^{i*}c_{ss}^{(j)}Z^i,$
$q^{(j)}=\sum_{s=1}^{R(j)}q_{ss}^{(j)}$ and
$Q=\sum_{j=1}^Nq^{(j)}=1_{D_1}.$ Note that
$Q=\sum_{i=0}^{k-1}(Z^i)^*qZ^i.$ Note that
\begin{align}\label{eMP5}
q_{ss}^{(j)}Z&=(\sum_{i=0}^{k-1}Z^{i*}c_{ss}^{(j)}Z^i)Z=Z\sum_{i=0}^k(Z^{i+1})^*c_{ss}^{(j)}Z^{i+1}\\
&=Z(\sum_{i=1}^{k-1}Z^{i*}c_{ss}^{(j)}Z^i+c_{ss}^{(j)})=Zq_{ss}^{(j)}.
\end{align}

It follows from \ref{GR} that $c_{11}^{(j)},$ $c_{11}^{(j)}Z^i$
and $c_{11}^{(j)}Z^kc_{11}^{(j)}$ generate a \SCA\, which is
isomorphic to $C(X_j)\otimes M_k$ for some compact subset
$X_j\subset S^1.$ Moreover, $q_{ss}^{(j)}Zq_{ss}^{(j)}$ is  in the
\SCA. Let $D$ be the \SCA\, generated by $D_1$ and
$c_{11}^{(j)}Z^kc_{11}^{(j)}.$ Then $D\cong\oplus_{j=1}^N
C(X_j)\otimes B_4\otimes M_k.$ It follows from  (\ref{eMP5}) that
$q^{(j)}$ and $Q$ commutes with $Z.$ Therefore $QZQ\in D.$ Thus,
by (\ref{eMps4}),
\begin{align}\label{e6}
\|Qu-uQ\| &\le  \|Qu-Qu_1\|+\|Qu_1-QZ\|+\|ZQ-u_1Q\|+\|u_1Q-uQ\|\\
&<2(\eta_3/32k +(3\pi/(k-2)+\eta_1/2))<\ep.\label{e6.5}
\end{align}
From $QZQ\in D,$ we also have
\begin{eqnarray}\label{e7}
QuQ\in_{\ep} D.
\end{eqnarray}
For $b\in {\cal F}_0,$ by (\ref{eMp1+}) and (\ref{eII}), we estimate  that
\begin{eqnarray}\nonumber
&&(Z^i)^*q(Z^i)b \approx_{2k\eta_1/2+k\eta_3/32k} (Z^i)^*qu_1^k
(W_1W_2\cdots W_i)(u^{k-i})^*b\\\nonumber
&=&(Z^i)^*qu_1^k(W_1W_2\cdots W_i)(u^{k-i})^*bu^{k-i}(W_1W_2\cdots
W_i)^* (u_1^k)^*[u^{k-i}(W_1W_2\cdots
W_i)^*(u_1^k)^*]^*\\\nonumber
\end{eqnarray}
Put $c_i=(u^{k-i})^*bu^{k-i}.$
Then
$c_i\in {\cal F}_1.$
Note that we have assumed that ${\cal F}_1\subset B_1.$
In particular, $p_1c_i=c_ip_1.$

Since $ (w_1w_2\cdots w_i){\cal F}_1(w_1w_2\cdots w_i)^*\subset
{\cal F}_3, $
 by (\ref{eMps3}),
\begin{eqnarray}\nonumber
&&(Z^i)^*qu_1^k(W_1W_2\cdots W_i)(u^{k-i})^*bu^{k-i}(W_1W_2\cdots
W_i)^* (u_1^k)^*[u^{k-i}(W_1W_2\cdots
W_i)^*(u_1^k)^*]^*\\\nonumber
&&=(Z^i)^*qu_1^k(W_1W_2\cdots W_i)c_i(W_1W_2\cdots W_i)^*
(u_1^k)^*[u^{k-i}(W_1W_2\cdots W_i)^*(u_1^k)^*]^*\\\nonumber
&&\approx_{\eta_0}(Z^i)^*qu_1^k(w_1w_2\cdots w_i)c_i(w_1w_2\cdots
w_i)^*u_1^k [u^{k-i}(W_1W_2\cdots W_i)^*(u_1^k)^*]^*\\\nonumber
&&\approx_{2k\eta_3/32k+\eta_1/32k+\eta_3/32k}
(Z^i)^*u_1^k(w_1w_2\cdots w_i)c_i(w_1w_2\cdots
w_i)^*u_1^kq[u^{k-i}(W_1W_2\cdots W_i)^*u_1^k]^*\\\nonumber
&&\approx_{\eta_0}(Z^i)^*u_1^k(W_1W_2\cdots
W_i)c_i(W_1W_2\cdots W_i)^* (u_1^{k})^*q[u^{k-i}(W_1W_2\cdots
W_i)^*(u_1^k)^*]^*\\\nonumber
&&
=(Z^i)^*u_1^k(W_1W_2\cdots
W_i)(u^{k-i})^*bu^{k-i}(W_1W_2\cdots W_i)^*
(u_1^{k})^*q[u^{k-i}(W_1W_2\cdots W_i)^*(u_1^k)^*]^*\\\nonumber
&&\approx_{2k\eta_3/32k}(Z^i)^*u_1^k(W_1W_2\cdots
W_i)(u_1^{k-i})^*bu_1^{k-i}(W_1W_2\cdots W_i)^*
(u_1^{k})^*q[u^{k-i}(W_1W_2\cdots W_i)^*(u_1^k)^*]^*\\\nonumber
&&\approx_{2k\eta_1/2}b(Z^i)^*qZ^i.
\end{eqnarray}
Note that ($k\ge 2+(4\cdot 256)\pi/\ep$)
\beq
(2k\eta_1/2+\eta_3/32)+\eta_0+\eta_1/32+\eta_1/32k+\eta_3/32k\\
+\eta_0+\eta_3/16+2k\eta_1/2<\ep/k^3.
\eneq
Hence
\begin{eqnarray}\label{eMp7}
\|(Z^i)^*qZ^ib-b(Z^i)^*qZ^i\|<\ep/k^3,\,\,\,k=0,1,...,k-1.
\end{eqnarray}
Therefore, for $b\in {\cal F}_0,$
\begin{eqnarray}\label{eMP8}
\|Qb-bQ\|<k(\ep/k^3)=\ep/k^2.
\end{eqnarray}
It follows from (\ref{e7}) and (\ref{e6.5}) that
\begin{eqnarray}\label{e9}
\|Qa-aQ\|<\ep\,\,\,\rforal\, a\in {\cal F}.
\end{eqnarray}

For any $b\in {\cal F}_0,$ a similar estimation above shows that
\begin{eqnarray}\label{eMP10}\nonumber
&&\hspace{-0.5in}\|qZ^ib(Z^i)^*q-qu_1^k(w_1w_2\cdots
w_i)(u^{k-i})^*bu^{k-i}(w_1w_2\cdots w_i)^*(u_1^k)^*q\|\\\nonumber
&&<2(2k\eta_1/2+k\eta_3/32k)=2k\eta_1+\eta_3/16.
\end{eqnarray}
However, by (\ref{eII}), (\ref{eMp1}), (b) and (\ref{eMp1}),
$$
qu_1^k(w_1w_2\cdots w_i)(u^{k-i})^*bu^{k-i}(w_1w_2\cdots
w_i)^*(u_1^k)^*q\in_{2\eta_3/32k+\eta_1/2+\eta_1/64+\eta_1/2}
B_4.
$$
It follows that, for $b\in {\cal F}_0,$
\begin{eqnarray}\label{eMp11}
(Z^i)^*qZ^ib(Z^i)^*qZ^i\in_{\ep/k^3}(Z^i)^*B_4Z^i,\,\,\,i=0,1,2,...,k-1.
\end{eqnarray}
Combing (\ref{eMp11}) with (\ref{eMp7}), we obtain that, for $b\in
{\cal F}_0,$
\begin{eqnarray}\label{eMP12}
QbQ\in_{k^2\ep/k^3}D_1\subset D.
\end{eqnarray}
Combing with (\ref{e7}), we obtain that
\begin{eqnarray}\label{eMP13}
QaQ\in_{\ep} D\rforal a\in {\cal F}.
\end{eqnarray}
We also compute that
\begin{eqnarray}\label{eMP14}\nonumber
[1-Q]&\le& [1-\sum_{i=1}^ke_i]+[1-p_1]\\
&\le& [r_2]+[r_1]\le [r_{00}]\le [b_0].
\end{eqnarray}
Combing (\ref{e9}), (\ref{eMP13}) and (\ref{eMP14}), by
\ref{one+0}, we conclude that $TR((A\rtimes_{\alpha}\Z)\otimes Q)=0.$  It follows from Theorem 3.6 of \cite{LN} that
$TR((A\rtimes_\af\Z)\otimes M)=0$ for any UHF-algebra $M$ of infinite type. On the other hand, by \cite{OP}, $A\rtimes_\af\Z$ has real real rank zero, stable rank one and weakly unperforated $K_0(A\rtimes_\af \Z).$ It follows from a classification result (Theorem 5.4 of \cite{LnNadv}) that $(A\rtimes_\af\Z)\otimes {\cal Z}$ is
isomorphic to  a unital simple \CA\, with tracial rank zero.
Since $\tau\circ \af^J=\tau$ for all $\tau\in T(A),$ the tracial cyclic Rokhlin property clearly implies the weak Rokhlin property.
By 4.9 of \cite{MS}, $A\rtimes_\af\Z$ is ${\cal Z}$ stable. It follows that $TR(A\rtimes_\af\Z)=0.$

\end{proof}

\begin{rem}\label{Reuni}
{\rm The proof presented above corrects 2.9 of \cite{LO}. If
$\alpha^k={\rm id}_A$ for some $k>1,$ however, the proof of 2.9 of
\cite{LO} works by replacing $2$ by $k$ (with $D\cong
e_1Ae_1\otimes M_k$). Since in this case the claim that $pwp\in D$
remains valid, conclusion of 2.9 of \cite{LO} holds.
One should also note that if  $\alpha^J$ is approximately inner for some integer $J>0,$ then
$[\alpha^J]=[{\rm id}_A]$ in $KL(A,A)$ An early version
of the proof of Theorem \ref{TM} contains an error which
was pointed us by Hiroki Matui. So it is appropriate to acknowledge
this at this point and  to choose their recent result (\cite{MS}) to simplify
the proof. More general result related to \ref{TM} will
be discussed elsewhere.
}
\end{rem}

\begin{prop}\label{Pdiv}
Let $A$ be a unital separable simple \CA\, with real rank zero,
stable rank one and weakly unperforated $K_0(A)$ with unique
tracial state. Let  $\alpha\in Aut(A)$ have the tracial cyclic
Rokhlin property.
Then there exists a subgroup $G\subset K_0(A)$
such that $\alpha^2_{*0}(g)=g$ for all $g\in G$ and $\rho_A(G)$ is
dense in $Aff(K_0(A)).$
\end{prop}

\begin{proof}
Let $G\subset K_0(A)$ such that
$\alpha^2_{*0}(g)=g$ for all $g\in G.$
Then $[1_A]\in G.$
Let $\tau\in T(A)$ be the
unique tracial state.
To show that $\rho_A(G)$ is dense in $\rho_A(K_0(A)),$ it suffices to show that, for any $n>0,$
there exists a nonzero projection
$e\in A$ such that $[\alpha^2(e)]=[e]$ and $\tau(e)<1/n.$

Therefore it suffices to show the following, for any $1>\ep>0$ and
any $0<\dt<{1-\ep\over{4}},$
if $p\in A$ is a nonzero projection such that
$$
\|\alpha^2(p)-p\|<\ep
$$
there is a nonzero projection $q\le p$ such that $p-q\not=0,$
$$
\|\alpha^2(q)-q\|<\ep+\dt \andeqn \|\alpha^2(p-q)-(p-q)\|<\ep+\dt.
$$
Because we must have that  $[p], [q], [p-q]\in G,$
that  either $\tau(p-q)\le {\tau(p)\over{2}}$ or $\tau(q)\le {\tau(p)\over{2}}$ and
that $\ep+\dt<1.$

There is at least one such $p$ (namely, $1_A$).

For any finite subset $(p\in ){\cal G}\subset A$ and $0<\eta<\min\{{\dt\over{256}}, {1-\ep\over{256}}, {\ep\over{256}}\},$
there are nonzero
mutually orthogonal projections $f_i, i=1,...,4$ such that

(i) $\|\alpha(f_i)-f_{i+1}\|<\eta,$ $i=1,2,...,4$ with  $f_5=f_1,$

(ii) $\|f_ia-af_i\|<\eta \rforal a\in {\cal G}$ and

(iii) $\tau(1-\sum_{i=1}^4f_i)<\eta.$

Put $e=f_1+f_3.$ Then
$$
\|\alpha^2(e)-e\|<2\eta.
$$

By taking sufficiently large ${\cal G}$ and
sufficiently small $\eta,$ by applyng \ref{Ad4}, there exist nonzero  projections $q, q_1\le p$
such that
$$
\|q-epe\|<2\eta\andeqn \|q_1-(1-e)p(1-e)\|<2\eta
$$
We compute that
\beq\nonumber
\|\alpha^2(q)-q\| &\le& \|\alpha^2(q)-\alpha^2(epe)\|+\|\alpha^2(epe)-epe\|\\\nonumber
&\le & 2\eta+\|\alpha^2(e)\alpha^2(p)\alpha^2(e)-epe\|<2\eta+2\eta+\eta+\ep<\ep +\dt/51<\ep+\dt \nonumber
\eneq
A similarly argument also shows that
\beq\nonumber
\|\alpha^2(p-q)-(p-q)\| <\ep+\dt.
\eneq

\end{proof}

\section{Tracial Rokhlin Property}

Let $A$ be a unital separable simple \CA\, with $TR(A)=0$ and let
$\alpha$ be an automorphism on $A.$  We now turn to the question
when $\alpha$ has the tracial cyclic Rokhlin property.

Let $X$ be a compact metric space, let $\sigma: X\to X$ be a
homeomorphism and let $\mu$ be a normalized $\sigma$-invariant
Borel measure. Recall that $\sigma$ has the Rokhlin property if,
for any $\ep>0$ and any $n\in {\mathbb N},$ there exists a Borel
set $E\subset X$ such that $E,\sigma(E),...,\sigma^n(E)$ are
mutually disjoint and
$$
\mu(X\setminus \cup_{i=0}^n\sigma^i(E))<\ep.
$$
From this, one may argue that the tracial Rokhlin property for
automorphism $\alpha$ above is  natural generalization of the
commutative case. In fact, from Theorem \ref{IT1}, it seems that
the tracial Rokhlin property occurs more often than one may first
thought and it appears that it is a rather natural phenomenon in
the context of automorphisms on simple \CA s.

Kishimoto originally studied approximately inner (but outer)
automorphisms regarding Problem {\bf P1}. It was recently proved
by N. C. Phillips that for a unital separable simple \CA\, with
tracial rank zero there is a dense $G_{\dt}$-set of approximately
inner automorphisms which satisfy the tracial Rokhlin property.
This shows that automorphisms with the tracial Rokhlin property
are abundant. But do they also have the tracial cyclic Rokhlin
property? It is proved (also follows from Theorem \ref{IIT2}
below) that if $\alpha$ is approximately inner as in the
Kishimoto's original case, tracial Rokhlin property implies the
tracial cyclic Rokhlin property.


\begin{lem}\label{IIL0}
Let $A$ be a unital separable simple \CA\, with real rank zero and
stable rank one and let $\alpha\in Aut(A).$ Suppose that there is
a subgroup $G\subset K_0(A)$ such that $\rho_A(G)$ is dense in
$\rho_A(K_0(A))$ and $(\alpha)_{*0}|_{G}={\rm id}_{G}.$ Then
$$
\tau(a)=\tau(\alpha(a))
$$
for all $a\in A.$

\end{lem}

\begin{proof}

Let $p\in A$ be a projection. Since $\rho_A(G)$ is dense in
$\rho_A(K_0(A))$ and $A$ is simple and has real rank zero and
stable rank one, there are projections $p_n\le p$ such that
$$
[p_n]\in G\andeqn \lim_{n\to\infty}\sup\{\tau(p-p_n):\tau\in
T(A)\}=0.
$$
Similarly, there are projections $e_n\le 1-(p-p_n)$ such that
$$
[e_n]\in G\andeqn
\lim_{n\to\infty}\sup\{\tau(1-(p-p_n)-e_n):\tau\in T(A)\}=0.
$$
It follows that
$$
\lim_{n\to\infty}\sup\{\tau(1-e_n):\tau\in T(A)\}=0.
$$
However, $[1-e_n]\in G,$ by the assumption,
$$
[\alpha(1-e_n)]=[1-e_n]\,\,\,{\rm in}\,\,\, K_0(A).
$$
Since $A$ has stable rank one,
 we have
$$
\tau(\alpha(1-e_n))=\tau(1-e_n)\to 0
$$
uniformly on $T(A).$ Since $p-p_n\le 1-e_n,$ we conclude that
$$
\tau(\alpha(p-p_n))\le \tau(\alpha(1-e_n))\to 0
$$
uniformly on $T(A).$ Thus, for any $\ep>0,$ there exists $N$ such
that
$$
\tau(p-p_n)<\ep/2\andeqn \tau(\alpha(p-p_n))<\ep/2
$$
for all $n\ge N$ and $\tau\in T(A).$ It follows that
\begin{eqnarray}
|\tau(\alpha(p))-\tau(p)|&\le &
|\tau(\alpha(p-p_n)|+|\tau(\alpha(p_n))-\tau(p_n)|+
|\tau(p-p_n)|\\
 &<&\ep/2+0+\ep/2=\ep.
 \end{eqnarray}
 Therefore, for any projection $p\in A,$
 $$
 \tau(\alpha(p))=\tau(p)
 $$
 for all $\tau\in T(A).$

Since $A$ has real rank zero, the above implies that
$\tau(\alpha(a))=\tau(a)$ for all $a\in A_{s.a}.$ The lemma then
follows.

\end{proof}

\begin{thm}\label{IIT1}
Let $A$ be a unital separable amenable simple \CA\, with $TR(A)=0$ which satisfies the UCT and let
$\alpha\in Aut(A).$ Suppose that $\alpha$ satisfies the tracial
Rokhlin property. If there is an integer $r>0$ such that
$\alpha^r_{*0}|_G={\rm id}_G$ for some subgroup $G\subset K_0(A)$
for which $\rho_A(G)$ is dense in $\rho_A(K_0(A)),$ then $\alpha$
satisfies the tracial cyclic Rokhlin property.
\end{thm}

Theorem \ref{IIT1} strengthens Theorem 3.14 of \cite{Lncd} slightly. This
is done by improving Lemma 6.3 of \cite{Lncd}. The rest of the
proof will be exactly the same  as that of  Theorem 3.14 of
\cite{Lncd} (which follows closely an idea of Kishimoto) but
applying Lemma \ref{IIL1} below. 


\begin{lem}\label{IILn}
Let $A$ be a unital simple separable \CA\, with $TR(A)=0$ and let 
$\af\in Aut(A)$ such that $(\af)_{*0}|_{G}={\rm id}_G$ for some 
subgroup $G\subset K_0(A)$ for which $\rho_A(G)={\rm Aff}(K_0(A)).$ Suppose that $\{p_j\}$ is a central sequence of projections such that $[p_j]\in G$ and define $\phi_j(a)=p_jap_j$ and $\psi_n(a)=\af(p_j)a\af(p_j),$ $j=1,2,....$ Then 
$\{\phi_j\}$ and $\{\psi_j\}$ are two sequentially asymptotic morphisms. Suppose also that there are finite-dimensional 
\SCA s $B_j$ and $C_j=\af(B_j)$ with $1_{B_j}=p_j$ and $1_{C_j}=\af(p_j)$ such that $[p_{j,i}]\in G$ for each minimal 
central projection $p_{j,i}$ of $B_j$ ($1\le i\le k(j)$) and there are sequentially asymptotic morphisms $\{\phi_j'\}$ and 
$\{\psi_j'\}$  such that
\beq\nonumber
&&\phi_j'(a)\subset B_j,\,\,\, \psi_j'(a)\subset C_j,\\\nonumber
&&\lim_{n\to\infty}\|\phi_j(a)-\phi_j'(a)\|=0\andeqn  \lim_{n\to\infty}\|\psi_j(a)-\psi_j'(a)\|=0
\eneq
for all $a\in A.$ Then, for any $\ep>0$ and  for any finite subset ${\cal G}\subset A$ and any finite subste of projections ${\cal P}_0\subset M_k(A)$ (for some $k\ge 1$) for which $[p]\in G$ for all $p\in {\cal P}_0,$ there exists an integer $J>0$ such that
$$
|\tau\circ \phi_j(a)-\tau\circ \psi_j(a)|<\ep/\tau(p_j)\rforal a\in {\cal G}
$$
and for all $\tau\in T(A),$ and, for all $j\ge J,$
$$
[\phi_j(p)]=[\psi_j(p)]\,\,\,{\rm in}\,\,\, K_0(A).
$$

\end{lem}

\begin{proof}
The proof is exactly the same as that of Lemma 6.2 of \cite{Lncd}. Instead of applying 
Lemma 6.1 of \cite{Lncd}, we apply \ref{IIL0}. 
\end{proof}

\begin{lem}\label{IIL1}
Let $A$ be a unital separable amenable simple  \CA\, with $TR(A)=0$
 which satisfies the UCT and let $\alpha\in
Aut(A)$ be such that $\alpha_{*0}|_G={\rm id}_{G}$ for some
subgroup $G$ of $K_0(A)$ for which $\rho_A(G)=\rho_A(K_0(A)).$
Suppose also that $\{p_j(l)\},$ $l=0,1,2,...,L,$ are central
sequences of projections in $A$ such that
$$
p_j(l)p_j(l')=0\,\,\,\text{if}\,\,\,l\not=l'\andeqn
\lim_{j\to\infty}\|p_j(l)-\alpha^l(p_j(0))\|=0,\,\,\, 1\le l\le L
$$
Then there exist central sequences of projections $\{q_j(l)\}$ and
 central sequences of partial isometries $\{u_j(l)\}$ such
that $q_j(l)\le p_j(l)$
$$
u_j(l)^*u_j(l)=q_j(0),\,\,\,u_j(l)u_j^*(l)=\alpha^l(q_j(0))
$$
for all large $j,$ and
$$
\lim_{j\to\infty}\|\alpha^l(q_j(0))-q_j(l)\|=0 \andeqn
\lim_{j\to\infty}\tau(p_j(l)-q_j(l))=0
$$
uniformly on $T(A).$
\end{lem}

\begin{proof}
 The proof is exactly the same as that of 6.3 of \cite{Lncd}.  On page 886 of that proof, 
 it uses the fact that $\rho(G)$ is dense in ${\rm Aff}(T(A))$ to produce projection $d_{n(j), t,s}.$
 This can be done with the current assumption which implies 
 that $\rho(G)$ is dense in ${\rm Aff}(T(A)).$ The proof of 6.3 of \cite{Lncd} also used Lemma 6.2 of \cite{Lncd}.
 At the end of p.887, one can apply \ref{IIL1} instead of 6.2 of \cite{Lncd}.

\end{proof}

Combining \ref{TM} with \ref{IIT1}, we have the following:

\begin{thm}\label{IIT2}
Let $A$ be a unital separable simple amenable \CA\, with $TR(A)=0$ which satisfies the UCT
and let $\alpha$ be an automorphism with the tracial Rokhlin
property. Suppose also that, for some integer $r>0,$
$\alpha^r_{*0}|_G={\rm id}|_G$ for some subgroup $G\subset K_0(A)$
for which $\rho_A(G)$ is dense in $\rho_A(K_0(A)).$ Then
$TR(A\rtimes_{\alpha}\Z)=0.$
\end{thm}

\begin{cor}\label{IIC1}
Let $A$ be a unital simple AH-algebra with slow dimension growth
and with real rank zero and let $\alpha\in Aut(A).$ Suppose that
$\alpha$ has the tracial Rokhlin property and
$[\alpha^r]=[{\rm id}_A]$ in $KL(A,A)$ for some integer $r\ge 1.$
 Then $A\rtimes_{\alpha}\Z$ is again an
AH-algebra  with slow dimension growth and with real rank zero.
\end{cor}

\begin{proof}
It is known that $TR(A)=0.$ By \ref{IIT2},
$TR(A\rtimes_{\alpha}\Z)=0.$ $A\rtimes_{\alpha}\Z$ also satisfies
the Universal Coefficient Theorem, by \cite{Lnduke},
$A\rtimes_{\alpha}\Z$ is an AH-algebra  with slow dimension growth
and with real rank zero.
\end{proof}

\vspace{0.2in}

\noindent
Department of Mathematics\\
East China Normal University\\
Shanghai, China\\
and (current)\\
Department of Mathematics\\
University of Oregon\\
Eugene, Oregon 97403, USA.

\end{document}